\documentclass[10pt,twoside,reqno]{amsart}
\usepackage{amssymb}
\usepackage{multirow}
\calclayout

\numberwithin{equation}{section}
\makeatletter

\renewcommand{\@secnumfont}{\bfseries}

\renewcommand{\section}{\@startsection{section}{1}%
  {0mm}{.7\linespacing\@plus\linespacing}{.5\linespacing}
  {\normalfont\bfseries\centering}}

\newcommand{\bibsection}{\@startsection{section}{1}%
  {0mm}{.7\linespacing\@plus\linespacing}{.5\linespacing}
  {\normalfont\scshape\centering}}

\renewcommand{\@biblabel}[1]{#1.}

\newtheorem{thm}{\bf Theorem}[section]

\begin{document}

\vspace{1.3cm}

\title {A note on Carlitz's type $q$-Changhee numbers and polynomials}

\author{Dmitry V. Dolgy}
\address{Hanrimwon, Kwangwoon University, Seoul 139-701, Republic
	of Korea}
\email{{d\_dol@mail.ru}}

\author{Gwan-Woo Jang}
\address{Department of Mathematics, Kwangwoon University, Seoul 139-701, Republic
	of Korea}
\email{gwjang@kw.ac.kr}

\author{Hyuck-In Kwon}
\address{Department of Mathematics, Kwangwoon University, Seoul 139-701, Republic
	of Korea}
\email{sura@kw.ac.kr}

\author{Taekyun Kim}
\address{Department of Mathematics, Kwangwoon University, Seoul 139-701, Republic of Korea}
\email{tkkim@kw.ac.kr}


\begin{abstract}
In this paper, we consider the Carlitz's type $q$-analogue of Changhee numbers and polynomials and we give some explicit formulae for these numbers and polynomials.
\end{abstract}
\maketitle
\bigskip
\medskip
\section{Introduction}
Let $p$ be an odd prime number. Throughout this paper, $\mathbb{Z}_p$, $\mathbb{Q}_p$ and $\mathbb{C}_p$ will denote the ring of $p$-adic integers, the field of $p$-adic numbers and the completion of the algebraic closure of $\mathbb{Q}_p$. The $p$-adic norm is normalized as $|p|_p = \tfrac{1}{p}$. Let $q$ be an indeterminate in $\mathbb{C}_p$ such that $|1-q|_p < p^{-\frac{1}{p-1}}$. The $q$-analogue of number $x$ is defined as $[x]_q = \frac{q^x-1}{q-1}$. As is well known, the Euler polynomials are defined by the generating function to be

\begin{equation}\begin{split}\label{01}
\frac{2}{e^t+1} e^{xt} = \sum_{n=0}^\infty E_n(x) \frac{t^n}{n!},\quad (\textnormal{see} \,\, [1-14]).
\end{split}\end{equation}
When $x=0$, $E_n=E_n(0)$, $(n \ge 0)$, are called the Euler numbers. In [1,2,3] L. Carlitz considered the $q$-analogue of Euler numbers which are given by the recurrence relation as follows:

\begin{equation*}\begin{split}
\mathcal{E}_{0,q}=1,\,\, q(q\mathcal{E}_q+1)^n + \mathcal{E}_{n,q} = \begin{cases}
[2]_q,&\text{if}\,\, n=0,\\
0,&\text{if}\,\,n>1. \end{cases}
\end{split}\end{equation*}
with the usual convention about replacing $\mathcal{E}_q^n$ by $\mathcal{E}_{n,q}$.

He also considered $q$-Euler polynomials which are defined by

\begin{equation}\begin{split}\label{02}
\mathcal{E}_{n,q}(x) = \sum_{l=0}^n {n \choose l} [x]_q^{n-l} q^{lx} \mathcal{E}_{l,q},\quad (\textnormal{see} \,\, [2,3]).
\end{split}\end{equation}
In [8,9,10], Kim defined the fermionic $p$-adic $q$-integral on $\mathbb{Z}_p$ as follow
s:

\begin{equation}\begin{split}\label{03}
I_{-q}(f) = \int_{\mathbb{Z}_p} f(x)   d\mu_{-q} (x) = \lim_{N \rightarrow \infty} \frac{1}{[p^N]_{-q}} \sum_{x=0}^{p^N-1} f(x) (-q)^x,
\end{split}\end{equation}
where $f(x)$ is continuous function on $\mathbb{Z}_p$ and $[x]_{-q} = \frac{1+q}{1-(-q)^x}.$

From \eqref{03}, He derived the following formula for the Carlitz's $q$-Euler polynomials:

\begin{equation}\begin{split}\label{04}
\int_{\mathbb{Z}_p} [x+y]_q^n   d\mu_{-q} (y) = \mathcal{E}_{n,q}(x),\,\,(n \geq 0),\quad (\textnormal{see} \,\, [7,10]).
\end{split}\end{equation}

When $x=0$, $\mathcal{E}_{n,q} = \int_{\mathbb{Z}_p} [x]_q^n d\mu_{-q} (x)$ are Carlitz's $q$-Euler n/umbers.

The Changhee polynomials are defined by the generating function to be

\begin{equation}\begin{split}\label{05}
\frac{2}{2+t}  (1+t)^x = \sum_{n=0}^\infty Ch_n(x) \frac{t^n}{n!},\quad (\textnormal{see} \,\, [5,6]).
\end{split}\end{equation}

Thus, by \eqref{05}, we get

\begin{equation}\begin{split}\label{06}
E_n(x) = \sum_{k=0}^n S_2(n,k) Ch_k(x),\,\, Ch_n(x) = \sum_{k=0}^n S_1(n,k) E_k(x),\,\,(n \geq 0),
\end{split}\end{equation}

Where $S_2(n,k)$ is Stirling number of the second kind and $S_1(n,k)$ is the Stirling number of the first kind. In [10], the higher-order Carlitz's $q$-Euler polynomials are written by the fermionic $p$-adic $q$-integral on $\mathbb{Z}_p$ as follows:

\begin{equation}\begin{split}\label{07}
\sum_{n=0}^\infty \mathcal{E}_{n,q}^{(r)}(x) \frac{t^n}{n!} =
\int_{\mathbb{Z}_p}\cdots\int_{\mathbb{Z}_p} e^{[x_1+\cdots+x_r+x]_q t}    d\mu_{-q} (x_1)\cdots    d\mu_{-q} (x_r),\,\,(n \geq 0).
\end{split}\end{equation}

In this paper, we consider the Carlitz's type $q$-Changhee polynomials and numbers and we give explicit formulas for these numbers and polynomials.

\section{Carlitz's type $q$-Changhee polynomials}

In this section, we assume that $t \in \mathbb{C}_p$ with $|t|_p < p^{-\frac{1}{p-1}}$. From \eqref{03} and \eqref{05}, we note that

\begin{equation}\begin{split}\label{08}
\int_{\mathbb{Z}_p} (1+t)^{x+y}d\mu_{-1}(y) = \frac{2}{2+t} (1+t)^x = \sum_{n=0}^\infty Ch_n(x) \frac{t^n}{n!},\quad (\textnormal{see} \,\, [4,5,6]).
\end{split}\end{equation}

Thus, by \eqref{08}, we get

\begin{equation}\begin{split}\label{09}
\int_{\mathbb{Z}_p} (x+y)_n    d\mu_{-1} (y) = Ch_n(x),\,\,(n \geq 0),
\end{split}\end{equation}

where $(x)_0=1$, $(x)_n = x(x-1)\cdots(x-n+1)$, $(n \geq 1)$.

In the viewpoint of \eqref{04}, we consider the Carlitz's type $q$-Changhee polynomials which are derived from the fermionic $p$-adic $q$-integral on $\mathbb{Z}_p$ as follows:

\begin{equation}\begin{split}\label{10}
\int_{\mathbb{Z}_p} (1+t)^{[x+y]_q}   d\mu_{-q} (y) = \sum_{n=0}^\infty Ch_{n,q}(x) \frac{t^n}{n!}.
\end{split}\end{equation}
Thus, by \eqref{10}, we get

\begin{equation}\begin{split}\label{11}
\sum_{n=0}^\infty Ch_{n,q}(x) \frac{t^n}{n!} &=
\sum_{k=0}^\infty \int_{\mathbb{Z}_p}  [x+y]_q^k   d\mu_{-q} (y) \frac{1}{k!} \Big( \log(1+t) \Big)^k \\
&=\sum_{k=0}^\infty \int_{\mathbb{Z}_p}  [x+y]_q^k   d\mu_{-q} (y) \sum_{n=k}^\infty S_1(n,k) \frac{t^n}{n!}\\
&= \sum_{n=0}^\infty \left( \sum_{k=0}^n \mathcal{E}_{k,q}(x) S_1(n,k) \right) \frac{t^n}{n!}.
\end{split}\end{equation}
Indeed,
\begin{equation*}\begin{split}
\sum_{k=0}^n S_1(n,k) \int_{\mathbb{Z}_p} [x+y]_q^k   d\mu_{-q} (y)& = \sum_{k=0}^n S_1(n,k) \frac{1}{(1-q)^k} \sum_{l=0}^k {k \choose l} q^{lx} (-1)^l \frac{[2]_q}{1+q^{l+1}}\\
&= [2]_q \sum_{k=0}^n \sum_{l=0}^k \frac{1}{(1-q)^k} {k \choose l} q^{lx} (-1)^l \frac{S_1(n,k)}{1+q^{l+1}}.
\end{split}\end{equation*}

Therefore, by \eqref{11}, we obtain the following theorem.

\begin{thm}
For $n \geq 0$, we have
\begin{equation*}\begin{split}
Ch_{n,q}(x) &= [2]_q \sum_{k=0}^n \sum_{l=0}^k \frac{1}{(1-q)^k} {k \choose l} q^{lx} (-1)^l \frac{S_1(n,k)}{1+q^{l+1}}\\
&= \sum_{k=0}^n  S_1(n,k)\mathcal{E}_{n,q}(x).
\end{split}\end{equation*}
\end{thm}

From \eqref{04}, we note that
\begin{equation}\begin{split}\label{12}
\sum_{n=0}^\infty \mathcal{E}_{n,q}(x) \frac{t^n}{n!} = \int_{\mathbb{Z}_p} e^{[x+y]_qt}   d\mu_{-q} (y).
\end{split}\end{equation}
By \eqref{12}, we get

\begin{equation}\begin{split}\label{13}
\sum_{k=0}^\infty Ch_{k,q}(x) \frac{1}{k!} \big(e^t-1\big)^k = \int_{\mathbb{Z}_p} e^{[x+y]_qt}   d\mu_{-q} (y) = \sum_{n=0}^\infty \mathcal{E}_{n,q}(x) \frac{t^n}{n!}
\end{split}\end{equation}

On the other hand,
\begin{equation}\begin{split}\label{14}
\sum_{k=0}^\infty Ch_{k,q}(x) \frac{1}{k!} \big(e^t-1\big)^k &= \sum_{k=0}^\infty Ch_{k,q}(x) \sum_{n=k}^\infty S_2(n,k) \frac{t^n}{n!}\\
&= \sum_{n=0}^\infty \left( \sum_{k=0}^n Ch_{k,q}(x) S_2(n,k) \right) \frac{t^n}{n!}.
\end{split}\end{equation}

Thus, by \eqref{13} and \eqref{14}, we get the following theorem.

\begin{thm}
For $n \geq 0$, we have
\begin{equation*}\begin{split}
\mathcal{E}_{n,q}(x) = \sum_{k=0}^n Ch_{k,q}(x) S_2(n,k).
\end{split}\end{equation*}
\end{thm}

From Theorem 1, we note that
\begin{equation}\begin{split}\label{15}
&\sum_{n=0}^\infty Ch_{n,q}(x) \frac{t^n}{n!}\\
 &= [2]_q \sum_{n=0}^\infty \left( \sum_{k=0}^n S_1(n,k) \frac{1}{(1-q)^k} \sum_{l=0}^k {k \choose l} q^{lx} (-1)^l \sum_{m=0}^\infty \Big(-q^{l+1} \Big)^m    \right) \frac{t^n}{n!}\\
 &=[2]_q \sum_{m=0}^\infty (-q)^m \sum_{n=0}^\infty \left( \sum_{k=0}^n S_1(n,k) [m+x]_q^k \right) \frac{t^n}{n!}\\
 &=[2]_q \sum_{m=0}^\infty (-q)^m \sum_{n=0}^\infty {[m+x]_q \choose n } t^n \\
 &=[2]_q \sum_{m=0}^\infty (-q)^m (1+t)^{[m+x]_q}.
\end{split}\end{equation}

Therefore, by \eqref{15}, we obtain the following theorem.

\begin{thm}
The generating function of the Carlitz's type $q$-Changhee polynomials is given by
\begin{equation*}\begin{split}
[2]_q \sum_{m=0}^\infty (-q)^m (1+t)^{[m+x]_q} = \sum_{n=0}^\infty Ch_{n,q}(x) \frac{t^n}{n!}.
\end{split}\end{equation*}
In particular, $x=0$, we have
\begin{equation*}\begin{split}
[2]_q \sum_{m=0}^\infty (-q)^m (1+t)^{[m]_q} = \sum_{n=0}^\infty Ch_{n,q} \frac{t^n}{n!}.
\end{split}\end{equation*}
\end{thm}

From \eqref{03}, we easily note that

\begin{equation}\begin{split}\label{16}
qI_{-q} (f_1) + I_{-q}(f) = [2]_q f(0), \,\,\text{where}\,\, f_1(x)=f(x+1).
\end{split}\end{equation}

Thus, by \eqref{16}, we get

\begin{equation}\begin{split}\label{17}
q\int_{\mathbb{Z}_p}(1+t)^{[x+1+y]_q}    d\mu_{-q} (y) + \int_{\mathbb{Z}_p} (1+t)^{[x+y]_q}   d\mu_{-q} (y) = [2]_q (1+t)^{[x]_q}.
\end{split}\end{equation}

By \eqref{10} and \eqref{17}, we get

\begin{equation}\begin{split}\label{18}
\sum_{n=0}^\infty \Big( qCh_{n,q}(x+1) + Ch_{n,q}(x) \Big) \frac{t^n}{n!} = [2]_q \sum_{n=0}^\infty \Big([x]_q \Big)_n \frac{t^n}{n!}.
\end{split}\end{equation}
Comparing the coefficients on the both sides of \eqref{18}, we get
\begin{equation}\begin{split}\label{19}
qCh_{n,q}(x+1) + Ch_{n,q}(x)  = [2]_q \Big([x]_q\Big)_n = [2]_q \sum_{l=0}^n S_1(n,l) [x]_q^l, \,\,(n \geq 0).
\end{split}\end{equation}
Therefore, we obtain the following theorem.

\begin{thm}
For $n \geq 0$, we have
\begin{equation*}\begin{split}
qCh_{n,q}(x+1) + Ch_{n,q}(x) = [2]_q \sum_{l=0}^n S_1(n,l) [x]_q^l.
\end{split}\end{equation*}
\end{thm}

From \eqref{10}, we have

\begin{equation}\begin{split}\label{20}
\sum_{n=0}^\infty \int_{\mathbb{Z}_p} { [x+y]_q \choose n}   d\mu_{-q} (y) t^n = \sum_{n=0}^\infty \frac{Ch_{n,q}(x)}{n!} t^n.
\end{split}\end{equation}

Thus, by \eqref{20}, we get

\begin{equation*}\begin{split}
\int_{\mathbb{Z}_p} { [x+y]_q \choose n}   d\mu_{-q} (y)  = \frac{Ch_{n,q}(x)}{n!} , \,\,(n \geq 0).
\end{split}\end{equation*}

Now, we observe that
\begin{equation}\begin{split}\label{21}
(1+t)^{[x+y]_q} = (1+t)^{[x]_q+q^x[y]_q}=(1+t)^{[x]_q} \cdot (1+t)^{q^x [y]_q}.
\end{split}\end{equation}
Thus, by \eqref{21}, we get

\begin{equation}\begin{split}\label{22}
\sum_{n=0}^\infty Ch_{n,q}(x) \frac{t^n}{n!} &= \int_{\mathbb{Z}_p} (1+t)^{[x+y]_q}   d\mu_{-q} (y)\\
& = \sum_{n=0}^\infty \sum_{k=0}^n S_1(n,k) \int_{\mathbb{Z}_p} \Big( [x]_q+q^x [y]_q \Big)^k   d\mu_{-q} (y) \frac{t^n}{n!}\\
&= \sum_{n=0}^\infty \left( \sum_{k=0}^n \sum_{l=0}^k {k \choose l} S_1(n,k) [x]_q^{k-l} q^{lx} \mathcal{E}_{l,q} \right) \frac{t^n}{n!}.
\end{split}\end{equation}

Therefore, by \eqref{22}, we obtain the following theorem.

\begin{thm}
For $n \geq 0$, we have
\begin{equation*}\begin{split}
Ch_{n,q}(x) =\sum_{k=0}^n \sum_{l=0}^k {k \choose l} S_1(n,k) [x]_q^{k-l} q^{lx} \mathcal{E}_{l,q}.
\end{split}\end{equation*}
\end{thm}

From \eqref{03}, we note that

\begin{equation}\begin{split}\label{23}
\int_{\mathbb{Z}_p} f(x)   d\mu_{-q} (x)
&= \lim_{N \rightarrow \infty} \frac{1}{[p^N]_{-q}} \sum_{x=0}^{p^N-1}f(x)(-q)^x\\
&= \lim_{N \rightarrow \infty} \frac{1}{[dp^N]_{-q}}  \sum_{x=0}^{dp^N-1}f(x)(-q)^x,
\end{split}\end{equation}
where $d \in \mathbb{N}$ with $d \equiv 1$ (mod 2). For $d \in \mathbb{N}$ with $d \equiv 1$ (mod 2), we have
\begin{equation}\begin{split}\label{24}
\int_{\mathbb{Z}_p} f(x)   d\mu_{-q} (y) = \lim_{N \rightarrow \infty} \frac{1}{[dp^N]_{-q}}\sum_{a=0}^{d-1} \sum_{x=0}^{p^N-1} f(a+dx)(-q)^{a+dx}.
\end{split}\end{equation}
By \eqref{24}, we get
\begin{equation}\begin{split}\label{25}
&\int_{\mathbb{Z}_p} (1+t)^{[x+y]_q}   d\mu_{-q} (y)\\
&= \frac{[2]_q}{[2]_{q^d}} \sum_{a=0}^{d-1} (-q)^a \int_{\mathbb{Z}_p} (1+t)^{[d]_q \big[\tfrac{a+x}{d}+y\big]_{q^d}}   d\mu_{-q^d} (y)\\
&= \frac{[2]_q}{[2]_{q^d}} \sum_{a=0}^{d-1} (-q)^a \sum_{k=0}^\infty [d]_q^k \int_{\mathbb{Z}_p} \big[ \tfrac{a+x}{d}+y \big]_{q^d}^k   d\mu_{-q} (y) \frac{1}{k!} \big( \log(1+t)\big)^k\\
&= \frac{[2]_q}{[2]_{q^d}} \sum_{a=0}^{d-1} (-q)^a \sum_{n=0}^\infty \left( \sum_{k=0}^n [d]_q^k \mathcal{E}_{k,q^d} \big( \tfrac{a+x}{d} \big) S_1(n,k) \right) \frac{t^n}{n!}\\
&=\sum_{n=0}^\infty \left( \frac{[2]_q}{[2]_{q^d}} \sum_{a=0}^{d-1}\sum_{k=0}^n (-q)^a [d]_q^k\mathcal{E}_{k,q^d} \big( \tfrac{a+x}{d} \big)S_1(n,k) \right) \frac{t^n}{n!}.
\end{split}\end{equation}

Therefore, by \eqref{10} and \eqref{25}, we obtain the following theorem.

\begin{thm}
For $n \geq 0$, we have
\begin{equation*}\begin{split}
Ch_{n,q}(x) = \frac{[2]_q}{[2]_{q^d}} \sum_{k=0}^n [d]_q^k\left(\sum_{a=0}^{d-1} (-q)^a \mathcal{E}_{k,q^d} \big( \tfrac{a+x}{d} \big)S_1(n,k) \right).
\end{split}\end{equation*}
\end{thm}

For $r \in \mathbb{N}$, the higher-order Carlitz's type $q$-Changhee polynomials are also given by the multivariate fermionic $p$-adic $q$-integral as follows:

\begin{equation}\begin{split}\label{26}
\int_{\mathbb{Z}_p}\cdots\int_{\mathbb{Z}_p} (1+t)^{[x_1+x_2+\cdots+x_r+x]_q}   d\mu_{-q} (x_1)\cdots    d\mu_{-q} (x_r) =  \sum_{n=0}^\infty Ch_{n,q}^{(r)}(x) \frac{t^n}{n!}.
\end{split}\end{equation}

Thus, we note that

\begin{equation}\begin{split}\label{27}
&\int_{\mathbb{Z}_p}\cdots\int_{\mathbb{Z}_p} (1+t)^{[x_1+x_2+\cdots+x_r+x]_q}   d\mu_{-q} (x_1)\cdots    d\mu_{-q} (x_r) \\
&= \sum_{k=0}^\infty \int_{\mathbb{Z}_p}\cdots\int_{\mathbb{Z}_p} [x_1+\cdots+x_r+x]_q^k   d\mu_{-q} (x_1)\cdots    d\mu_{-q} (x_r) \frac{1}{k!} \Big( \log(1+t)\Big)^k\\
&= \sum_{n=0}^\infty \left( \sum_{k=0}^n S_1(n,k) \int_{\mathbb{Z}_p}\cdots\int_{\mathbb{Z}_p} [x_1+\cdots+x_r+x]_q^k   d\mu_{-q} (x_1)\cdots    d\mu_{-q} (x_r) \right) \frac{t^n}{n!}\\
&= \sum_{n=0}^\infty \left( \sum_{k=0}^n S_1(n,k) \mathcal{E}_{k,q}^{(r)}(x) \right) \frac{t^n}{n!}.
\end{split}\end{equation}
By \eqref{26} and \eqref{27}, we get

\begin{equation}\begin{split}\label{28}
Ch_{n,q}^{(r)}(x) = \sum_{k=0}^n S_1(n,k) \mathcal{E}_{k,q}^{(r)}(x).
\end{split}\end{equation}
When $x=0$, $Ch_{n,q}^{(r)}=Ch_{n,q}^{(r)}(0)$ are called the Carlitz's type $q$-Changhee numbers.

By \eqref{07} and \eqref{26}, we get
\begin{equation}\begin{split}\label{29}
\sum_{k=0}^\infty Ch_{k,q}^{(r)}(x) \frac{1}{k!} \Big(e^t-1\Big)^k &= \int_{\mathbb{Z}_p}\cdots\int_{\mathbb{Z}_p} e^{[x_1+x_2+\cdots+x_r+x]_q t} d\mu_{-q} (x_1)\cdots    d\mu_{-q} (x_r) \\
&= \sum_{n=0}^\infty \mathcal{E}_{n,q}^{(r)}(x).
\end{split}\end{equation}

On the other hand,

\begin{equation}\begin{split}\label{30}
\sum_{k=0}^\infty Ch_{k,q}^{(r)} (x) \frac{1}{k!} \Big(e^t-1\Big)^k &=\sum_{k=0}^\infty Ch_{k,q}^{(r)} (x)  \sum_{n=k}^\infty S_2(n,k) \frac{t^n}{n!}\\
&= \sum_{n=0}^\infty \left( \sum_{k=0}^n Ch_{k,q}^{(r)}(x) S_2(n,k) \right) \frac{t^n}{n!}.
\end{split}\end{equation}

Comparing the coefficients on the both sides of \eqref{29} and \eqref{30}, we have

\begin{equation}\begin{split}\label{31}
\mathcal{E}_{n,q}^{(r)}(x)= \sum_{k=0}^n Ch_{k,q}^{(r)}(x) S_2(n,k) .
\end{split}\end{equation}

\end{document}